\documentclass{aip-cp}
\usepackage[numbers]{natbib}
\usepackage{graphicx}
\usepackage{color, multirow,amsfonts}
\usepackage{url}
\usepackage{rotating}
\begin{document}
\title{Distribution sensitive estimators of the index of regular variation based on ratios of order statistics}
\author[aff1]{Pavlina K. Jordanova\corref{cor1}}
\author[aff2,aff3,aff4]{Milan Stehl\'\i k}
\eaddress{mlnstehlik@gmail.com}

\affil[aff1]{Faculty of Mathematics and Informatics, Konstantin Preslavsky University of Shumen, \\115 "Universitetska" str., 9712 Shumen, Bulgaria.}
\affil[aff2]{Department of Statistics and Actuarial Science, The University of Iowa, Iowa City, Iowa, USA.}
\affil[aff3]{Institute of Statistics, Universidad de Valpara\'iso, Valpara\'iso, Chile.}
\affil[aff4]{Department of Applied Statistics, Johannes Kepler University, Altenbergerstrasse 69, 4040 Linz, Austria.}

\corresp[cor1]{Corresponding author: pavlina\_kj@abv.bg}

\maketitle

\begin{abstract}
Ratios of central order statistics seem to be very useful for estimating the tail of the distributions and therefore, quantiles outside the range of the data. In 1995 Isabel Fraga Alves investigated the rate of convergence of three semi-parametric estimators of the parameter of the tail index in case when the cumulative distribution function of the observed random variable belongs to the max-domain of attraction of a fixed Generalized Extreme Value Distribution. They are based on ratios of specific linear transformations of two extreme order statistics. In 2019 we considered Pareto case and found two very simple and unbiased estimators of the index of regular variation. Then, using the central order statistics we showed that these estimators have many good properties. Then, we observed that although the assumptions are different, one of them is equivalent to one of Alves's estimators. Using central order statistics we proved unbiasedness, asymptotic consistency, asymptotic normality and asymptotic efficiency.  Here we use again central order statistics and a parametric approach and obtain distribution sensitive estimators of the index of regular variation in some particular cases. Then, we find conditions which guarantee that these estimators are unbiased, consistent and asymptotically normal. The results are depicted via simulation study.
\end{abstract}

\section{INTRODUCTION AND PRELIMINARIES}

Let us assume that $X_1, X_2, ..., X_n$ are $n$ independent observations on a random variable (r.v.) $X$ with cumulative distribution function (c.d.f.) $F_X(x) = \mathbb{P}(X \leq x)$ with regularly varying right tail. More precisely, we suppose that for some $\alpha > 0$,
$$\lim_{t \to \infty} \frac{\mathbb{P}(X > tx)}{\mathbb{P}(X > t)} =  \lim_{t \to \infty} \frac{\bar{F}_X(tx)}{\bar{F}_X(t)} = x^{-\alpha}, \quad \forall x > 0,$$
where $\bar{F}_X(x) = 1 - F_X(x)$.  Briefly we will denote these limit relations in this way $\bar{F}_X \in RV_{-\alpha}$.
The distributions of order statistics  $X_{(1,n)} \leq X_{(2,n)}  \leq ... \leq X_{(n,n)}$ are very well investigated in the scientific literature. One can see for example  Wilks (1948)\cite{Wilks1948}, Renyi (1953)\cite{Renyi1953}, Arnold (1992-2015)\cite{arnold1992first,Arnold2015}, or Nevzorov (2001)\cite{nevzorov2001records}.

The task for estimation of the index of regular variation $\alpha > 0$ has received much attention during the last years. If we use a non-parametric approach the well-known results seems to be quite rough and, therefore, useful only in cases of large samples. Here, we concentrate our study on a parametric approach. It turns out that again in very particular cases, for example in Pareto case, it already has a satisfactory solution. In 1995 Isabel Fraga Alves \cite{alves1995estimation} investigated the rate of convergence of three semi-parametric estimators of the tail index in case when the c.d.f. of the observed r.v. belongs to the max-domain of attraction of a fixed Generalized Extreme Value Distribution. One of these estimators is $Q_{i,s}^*$, defined in (\ref{QPar}). She works mainly with extreme order statistics. Here we consider central order statistics and the following estimators
\begin{eqnarray}
\label{QPar}  Q_{i,s}  &:=& \frac{\log\frac{X_{(is, n)}}{X_{(i, n)}}}{H_{is-1} - H_{i-1}}, \,\, Q_{i,s}^* := \frac{\log\frac{X_{(is, n)}}{X_{(i, n)}}}{\log(s)}, \\
\label{QstarLL}  Q_{i,s}^{LL}  &:=& \frac{Q_{i,s}}{2},\quad \quad \quad \quad  Q_{i,s}^{LL*} = \frac{Q_{i,s}^*}{2},\\
\label{QstarFr}  Q_{i,s}^{Fr*} &:=& -\frac{\log\frac{X_{(is, n)}}{X_{(i, n)}}}{\log\left[1-\frac{\log(s)}{\log(s+1)}\right]}, \\
\label{QstarHH}  Q_{i,s}^{HH*} &:=& \frac{\log\frac{X_{(is, n)}}{X_{(i, n)}} + \log\left[1 - \frac{\log(s)}{\log(s + 1)}\right]}{\log(s)},
\end{eqnarray}
where $s = 2, 3, ...$, $n = i(s+1)-1$, and $i = 2, 3, ...$. The fact that these estimators are functions only of ratios of order statistics entails the invariance of these estimators with respect to a deterministic scale change of the sample. Therefore, without lost of generality,  everywhere in this work we assume that the scale parameters $\sigma$ in the considered distributions  are equal to $1$. 

In 2019 Jordanova and Stehlik \cite{JordanovaStehlikAmitans2019}, assumed that the observed r.v. is Pareto distributed. Then, they have used the well-known formulae for the mean and the variance of logarithmic differences of the corresponding order statistics, and proved that $Q_{i,s}$, and $Q_{i,s}^*$   estimators are unbiased (the second one only asymptotically), consistent, asymptotically efficient, and asymptotic normal. Due to the fact that we fix the exact probability type of the observed r.v. they do not impose separately the second order regularly varying condition defined in Geluk et al. (1997) \cite{geluk1997second}. For Pareto distribution it is automatically satisfied. Here we follow the same approach for different probability types. Its main advantage is that it is very flexible and provides an useful accuracy given mid-range and small samples. Log-Logistic, Fr$\acute{e}$chet and Hill-horror cases are partially investigated. The conducted simulation study depicts the quality of the results.

Further on we denote by
$$H_{n, m} =  1 + \frac{1}{2^m} + \frac{1}{3^m} + ... + \frac{1}{(n-1)^m} + \frac{1}{n^m}, \quad n = 1, 2, ..., $$
the $n$-th Generalised harmonic number of power $m = 1, 2, ...$ and $H_n := H_{n, 1}$, denotes the $n$-th harmonic number.

Along the paper we denote by  $F^\leftarrow(p) = \inf \{x \in \mathbb{R}: F(x) \geq p\}$, the theoretical left-continuous version of the
quantile function of a c.d.f. $F$, for $p \in (0, 1]$ and by assumption $F^\leftarrow(0) := \sup \{x \in R: F(x) = 0\}$, $\sup \emptyset = -\infty$. The following definition of empirical quantile function $F_n^\leftarrow(0):= X_{(1, n)}$,
\begin{equation}\label{empquant3}
F_n^\leftarrow(p) = \left\{\begin{array}{ccc}
                             X_{([np+1], n)} & , & np \not\in  \mathbb{N}, \\
                             X_{(np, n)} & , & np \in \mathbb{N},
                           \end{array}
\right. =
\left\{\begin{array}{ccc}
                             X_{(\lceil np \rceil, n)} & , & np \not\in  \mathbb{N}, \\
                             X_{(np, n)} & , & np \in \mathbb{N},
                           \end{array}
\right. \quad p \in (0; 1],
\end{equation}
where $[a]$ means the integer part of $a$ and $\lceil a \rceil$ is for the ceiling of $a$, i.e. the least integer greater than or equal to $a$, could be seen e.g. in Serfling (2009) \cite{serfling2009approximation}. It is equivalent to the Definition 1, in Hyndman and Fan (1996) \cite{hyndman1996sample} and is implemented in function $quantile$ in software R \cite{R}, with parameter $Type = 1$.

\section{LOG-LOGISTIC CASE}

In this section we consider a sample of $n$ independent observations on a r.v. $X$ with c.d.f.
\begin{equation}\label{LogLogisticCDF}F_X(x) = \frac{1}{1+x^{-\alpha}},\quad  x > 0. \end{equation}
 Briefly $X \in Log-Logistic(\alpha; 0, 1)$.

Balakrishnan et al. (1987) \cite {balakrishnan1987best} found the best linear unbiased estimators of its location and
scale parameters. However, in their work the estimator of the shape parameter is missing. Recently Ahsanullah and Alzaatreh (2018) \cite{ahsanullah2018parameter} have considered again only estimation procedure of these two parameters and the ideas about the estimation of the shape parameter are reduced only to Hill's estimator (1975) \cite{hill1975simple}.

Here we propose two new estimators $ Q_{i,s}^{LL}$  and $ Q_{i,s}^{LL*}$ and investigate their properties. The idea comes from  the following considerations. By formula (\ref{LogLogisticCDF}) one can see that $X^\alpha \in Log-Logistic(1; 0, 1)$, and
therefore, the r.v. $\alpha\log(X)$ is standard Logistic distributed, i.e. its location parameter is 0, its scale and shape parameters are equal to 1.
The last conclusion allows us to reduce the task for estimation of the parameter $\alpha$, to the one of investigation of the properties of order statistics in the Logistic case which are already very well investigated in the scientific literature. In 1963 Birnbaum and Dudman\cite{birnbaum1963logistic} found their moment and  cumulant generating functions. Then, by using their derivatives the authors expressed the mean and the variance of Logistic order statistics via polygamma function correspondingly of the first and the second order. Gupta and Balakrishnan (1991) \cite{GuptaBala1991logistic} made a step further on and found a series representation of the joint moments of these order statistics. Their closed form solution seems to be still an open problem.
The following result is an immediate corollary of Theorem 4.6,a), in Jordanova (2020) \cite{Jordanova2020Monograph}.

{\bf Theorem 1.}\label{thmTh1LLcase} Let us fix $s = 2, 3, ...$ and consider a sample of $i(s+1)-1$, $i \in \mathbb{N}$ independent observations on a r.v. $X \in$ $Log-Logistic(\alpha; 0, 1)$, $\alpha > 0$. Then, for $x > 0$ the probability density function (p.d.f.) of $\alpha Q_{i,s}^{LL*}$ is
\begin{eqnarray}
  \label{withouthypergeometric}
   \nonumber f_{\alpha Q_{i,s}^{LL*}}(x) &=& 2\log(s)\frac{[i(s+1)-1]!s^{-2ix}(1-s^{-2x})^{i(s-1)-1}}{[(i-1)!]^2[i(s-1)-1]!}\int_0^{\infty}\frac{z^{is-1}}{(1+z)^{is}(1+zs^{-2x})^{is}}dz, \\
   \nonumber    &=& 2\log(s)\frac{[i(s+1)-1]![(is-1)!]^2s^{2xi(s-1)}(1-s^{-2x})^{i(s-1)-1}}{[i(s-1)-1]![(i-1)!]^2(2is-1)!} {}_2F_1(is,is;2is;1-s^{2x})
     \end{eqnarray}
  where ${}_2F_1(a,b;c;y) = \sum_{n=0}^\infty \frac{(a)_n(b)_n}{(c)_n}\frac{y^n}{n!} = \int_0^1 \frac{x^{b-1}(1-x)^{c-b-1}}{(1-yx)^aB(b,c-b)} dx$, $(q)_n = q(q+1)...(q+n-1)$, $(q)_0 = 1$, $c > b > 0$, $y < 1$, is the Gauss hypergeometric function (the Euler type integral) and $f_{\alpha Q_{i,s}^{LL*}}(x) =  0$, otherwise.

\medskip

 When compute the quantile function of the r.v. $X$ with c.d.f. (\ref{LogLogisticCDF}) we observe that for all $s \in \mathbb{N}$,
$$\log\left(\frac{F_X^\leftarrow(\frac{s}{s+1})}{F_X^\leftarrow(\frac{1}{s+1})}\right) = \frac{2}{\alpha}\log(s).$$
In Jordanova and Stehl\'\i k (2019), in general (not only for Log-Logistic) case we have shown that,
\begin{equation}\label{limitinprobability}
\log\left(\frac{X_{(is, i(s+1)-1)}}{X_{(i, i(s+1)-1)}}\right) {\mathop{\to}\limits_{i \to \infty}^{\mathbb{P}}} \log\left(\frac{F_X^\leftarrow(\frac{s}{s+1})}{F_X^\leftarrow(\frac{1}{s+1})}\right) , \quad  s \in \mathbb{N}
\end{equation}
Therefore, in this case
\begin{equation}\label{LLSrtongconsistency}
\log\left(\frac{X_{(is, i(s+1)-1)}}{X_{(i, i(s+1)-1)}}\right) {\mathop{\to}\limits_{i \to \infty}^{\mathbb{P}}} \frac{2}{\alpha}\log(s), \quad  s \in \mathbb{N}
\end{equation}
and when we normalize $\log\left(\frac{X_{(is, i(s+1)-1)}}{X_{(i, i(s+1)-1)}}\right)$ with  $2\log(s)$ we obtain a consistent estimator for $\frac{1}{\alpha}$. The last one is exactly  $Q_{i,s}^{LL*}$ estimator, defined in (\ref{QstarFr}). In order to obtain $Q_{i,s}^{LL}$ estimator we normalize $\log\left(\frac{X_{(is, i(s+1)-1)}}{X_{(i, i(s+1)-1)}}\right)$ with its expectation. The proof of the following theorem could be found in Jordanova (2020) \cite{Jordanova2020Monograph}. The approach is analogous to the one used in Jordanova and Stehl\'\i k (2019) \cite{JordanovaStehlikAmitans2019} for estimation of the parameter $\alpha$ in Pareto case.

{\bf Theorem 2.}\label{LLstrongconsistency} Let us fix $s = 2, 3, ...$ and consider a sample of $i(s+1)-1$, $i \in \mathbb{N}$ independent observations on a r.v. $X \in$ $Log-Logistic(\alpha; 0, 1)$, $\alpha > 0$.
\begin{description}
  \item[i)]  For all $i = 2, 3, ...$, $\mathbb{E}[Q_{i,s}^{LL}] = \frac{1}{\alpha}$.
  \item[ii)] $Q_{i,s}^{LL}  {\mathop{\to}\limits_{i \to \infty}^{a.s.}} \frac{1}{\alpha}$.
  \item[iii)]   $\mathbb{E} (Q_{i,s}^{LL*}) = \frac{H_{is-1} - H_{i-1}}{\alpha \,\,\log(s)} {\mathop{\to}\limits_{i \to \infty}^{}} \frac{1}{\alpha}.$
  \item[iv)]  $Q_{i,s}^{LL*}  {\mathop{\to}\limits_{i \to \infty}^{a.s.}} \frac{1}{\alpha}$.
  \item[v)] $2\sqrt{i(s+1) - 1}(H_{is-1} - H_{i-1})\left[\alpha Q_{i,s}^{LL}  - \frac{\log(s)}{H_{is-1} - H_{i-1}}\right] {\mathop{\to}\limits_{i \to \infty}^{d}} \eta,$   where $\eta \in N\left(0, \frac{2(s+1)^2(s-1)}{s^2}\right).$
  \item[vi)] $\sqrt{i(s+1) - 1}\left[\alpha Q_{i,s}^{LL*}  - 1\right] {\mathop{\to}\limits_{i \to \infty}^{d}} \eta,$
  where $\eta \in N\left(0, \frac{(s+1)^2(s-1)}{2s^2[\log(s)]^2}\right).$
  \end{description}

{\bf{Remarks.}}  The first statement in Theorem 2 means that for all $i \in \mathbb{N}$ and $s = 2, 3, ...$, $Q_{i,s}^{LL}$ are unbiased estimators for $\frac{1}{\alpha}$.  The limit relations ii) and iv) say that both $Q_{i,s}^{LL}$ and $Q_{i,s}^{LL*}$  are a strongly consistent when $i \to \infty$. The third one states that $Q_{i,s}^{LL*}$ is asymptotically unbiased estimator for $\frac{1}{\alpha}$. Points v) and vi) clarify that both $Q_{i,s}^{LL}$ and $Q_{i,s}^{LL*}$ are asymptotically normal when $i$ increases unboundedly.

The limit relation  vi) in Theorem 2 allows us to obtain large sample confidence intervals for $\alpha$. It means that
$$\frac{s\log(s)}{s+1}\sqrt{\frac{i(s+1) - 1}{\frac{s-1}{2}}}\left(\alpha Q_{i,s}^{LL*} - 1\right) {\mathop{\to}\limits_{i \to \infty}^{d}} \theta, \quad \theta \in N\left(0, 1\right).$$
Thus, if we choose confidence level $(1-\alpha_0)100\%$, where $\alpha_0 \in (0, 1)$ and if we denote by $z_{1-\frac{\alpha_0}{2}}$ the $1-\frac{\alpha_0}{2}$ quantile of the Standard Normal distribution, then
$$\lim_{i \to \infty}\mathbb{P}\left[-z_{1-\frac{\alpha_0}{2}} \leq \frac{s\log(s)}{s+1}\sqrt{\frac{2[i(s+1) - 1]}{s-1}}\left(\alpha Q_{i,s}^{LL*} - 1\right)
  \leq z_{1-\frac{\alpha_0}{2}}\right] = 1-\alpha_0.$$
Therefore,
$$\lim_{i \to \infty}\mathbb{P}\left[\frac{1}{Q_{i,s}^{LL*}}-\frac{z_{1-\frac{\alpha_0}{2}}(s+1)}{Q_{i,s}^{LL*} s\log(s)} \sqrt{\frac{s-1}{2[i(s+1) - 1]}}\leq \alpha  \leq \frac{1}{Q_{i,s}^{LL*}} + \frac{z_{1-\frac{\alpha_0}{2}}(s+1)}{Q_{i,s}^{LL*} s\log(s)} \sqrt{\frac{s-1}{2[i(s+1) - 1]}}\right] = 1-\alpha_0.$$

Now Jordanova (2020) \cite{Jordanova2020Monograph} uses the definition of $Q_{i,s}^{LL*}$ and concludes that for any fixed $s = 2, 3, ...$, and $i \in \mathbb{N}$ large enough, the corresponding $1-\alpha_0$-confidence intervals (c.is.) for $\alpha$ are:
\begin{equation}\label{ciLL}
\left[\frac{2\log(s)}{\log\frac{X_{(is, i(s+1)-1)}}{X_{(i, i(s+1)-1)}}} -\frac{2z_{1-\frac{\alpha_0}{2}}(s+1)}{s\log\frac{X_{(is, i(s+1)-1)}}{X_{(i, i(s+1)-1)}}} \sqrt{\frac{s-1}{2[i(s+1) - 1]}}; \frac{2\log(s)}{\log\frac{X_{(is, i(s+1)-1)}}{X_{(i, i(s+1)-1)}}}  +\frac{2z_{1-\frac{\alpha_0}{2}}(s+1)}{s\log\frac{X_{(is, i(s+1)-1)}}{X_{(i, i(s+1)-1)}}}\sqrt{\frac{s-1}{s[i(s+1) - 1]}}\right].
\end{equation}

The task for estimation of the quantiles outside the range of the data seems to be more difficult. It is easy to see that if $p < \frac{1}{n}$, then by using the definition (\ref{empquant3}) we will obtain one and the same estimator for $F_X^\leftarrow(1 - p)$ and it is $X_{(n,n)}$. In order to improve it we use formula (\ref{LogLogisticCDF}) and more precisely its inverse function $F_X^\leftarrow(p) = \left(\frac{p}{1 - p}\right)^{1/\alpha}$, and obtain the following two estimators
\begin{equation}\label{QuantileestimatorLL}
\hat{F}_{i,s,LL}^\leftarrow(p) = \left(\frac{p}{1 - p}\right)^{Q_{i,s}^{LL}}, \quad \quad \left(\hat{F}_{i,s,LL}^\leftarrow\right)^*(p) = \left(\frac{p}{1 - p}\right)^{Q_{i,s}^{LL*}}.
\end{equation}
Although they are not asymptotically normal, the following theorem explains why they are better than $X_{(n,n)}$ for estimation for $1-p$-th quantile, $F_X^\leftarrow(1 - p)$, $p < \frac{1}{n}$, which is outside the range of the data. Its statements follow by formulae (\ref{QuantileestimatorLL}), Theorem 2, ii) and iv), continuity of the function $g(x) = \left(\frac{p}{1 - p}\right)^x$ and Continuous mappings theorem.

{\bf Theorem 3.}\label{QuantileEstimatorsLLproperties} Let us fix $s = 2, 3, ...$ and consider a sample of $n = i(s+1)-1$, $i \in \mathbb{N}$ independent observations on a r.v. $X \in$ $Log-Logistic(\alpha; 0, 1)$, $\alpha > 0$. Then
\begin{description}
  \item[i)]   $\hat{F}_{i,s,LL}^\leftarrow(p) {\mathop{\to}\limits_{i \to \infty}^{a.s.}} F_X^\leftarrow\left(p\right)$.
  \item[ii)]  $\left(\hat{F}_{i,s,LL}^\leftarrow\right)^*(p)  {\mathop{\to}\limits_{i \to \infty}^{a.s.}} F_X^\leftarrow\left(p\right)$.
     \end{description}

{\it Simulation study}

\begin{figure}
\begin{minipage}[t]{0.5\linewidth}
    \includegraphics[scale=.58]{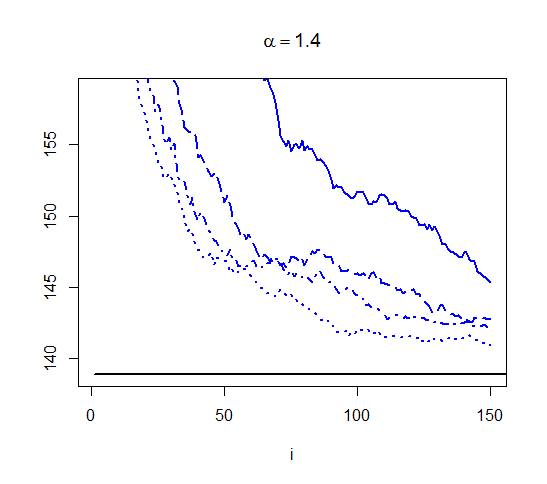}\vspace{-0.3cm}
\end{minipage}
\begin{minipage}[t]{0.49\linewidth}
    \includegraphics[scale=.58]{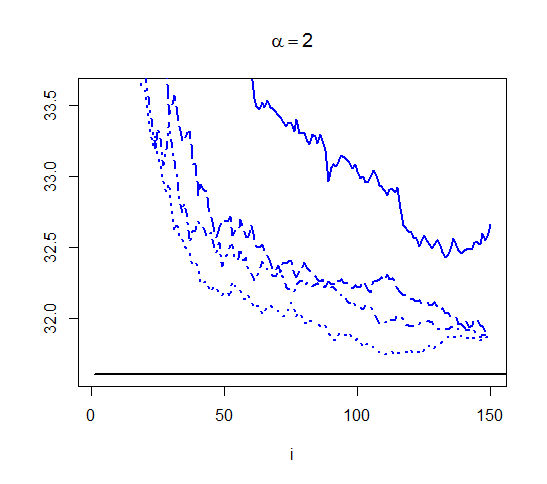}\vspace{-0.3cm}
\end{minipage}
\caption{\small Dependence of $\left(\hat{F}_{i,s,LL}^\leftarrow\right)^*(p)$ on $i$ in Log-Logistic case.\label{fig:LL0AMiTaNSquantiles}}
\end{figure}

The rate of convergence of $Q_{i,s}^{LL*}$ to $\frac{1}{\alpha}$ is quite good and it is partially depicted in Jordanova (2020) \cite{Jordanova2020Monograph}. For $\alpha < 1$ and $p = \frac{999}{1000}$ given small samples\footnote{Here the sample size $n < 899$ and $p = 999/1000$, which mean that $\left(\hat{F}_{i,s,LL}^\leftarrow\right)^*(p)$ is an estimator of a quantile outside the range of the data.} $\left(\hat{F}_{i,s,LL}^\leftarrow\right)^*(p)$ estimators are quite rough, therefore, here we have skipped their plots. The plots in Figure \ref{fig:LL0AMiTaNSquantiles} depict the rate of convergence of $\left(\hat{F}_{i,s,LL}^\leftarrow\right)^*(p)$ to $F_X^\leftarrow\left(p\right)$ for different $s = 2$ (solid lines), $s =  3$ (dashed lines), $s =  4$ (dash-dot lines), $s =  5$ (dotted lines), $\alpha = 1.4$ or $\alpha = 2$. The estimated value of $F_X^\leftarrow\left(p\right)$ is presented via a straight solid line. In order to plot these figures for different but fixed values of $s$ and $\alpha$, we have simulated $1000$ samples of $n = 150(s + 1) - 1$ independent observations on a r.v. $X \in Log-Logistic(\alpha; 0, 1)$. Then, for any fixed sample and for any fixed $i = 1, 2, ..., 150$
we have computed $\left(\hat{F}_{i,s,LL}^\leftarrow\right)^*(p)$, $p = \frac{999}{1000}$, based on formula (\ref{QuantileestimatorLL}). Finally, we have averaged $\left(\hat{F}_{i,s,LL}^\leftarrow\right)^*(p)$ over these $1000$ samples and we have plotted these averages as a function of $i = 10, 11, ..., 150$. We have two parameters which govern the sample size $n = i(s+1)-1$. These are $i$ and $s$. We observe that when the sample size increases the estimators $\left(\hat{F}_{i,s,LL}^\leftarrow\right)^*(p)$ get closer to the estimated value $F_X^\leftarrow\left(p\right)$.   These harmonize with our results in  Theorem 3, ii) which says that $\left(\hat{F}_{i,s,LL}^\leftarrow\right)^*(p)$ is a strongly consistent estimator for $F_X^\leftarrow\left(p\right)$.

\section{FR$\acute{E}$CHET CASE}
Let us now assume that the observed r.v. has c.d.f. $F_X(x) = 0$ for $x \leq 0$ and if $x > 0$, $F_X(x) = \exp\{-x^{-\alpha}\}$.
Briefly we will denote this by $X \in Fr\acute{e}chet(\alpha, 0, 1)$. In this case, $X^\alpha \in Fr\acute{e}chet$ $(1,$ $0$,$ 1)$ and $\alpha\log(X) \in Gumbel(0, 1)$. If we try to obtain Maximum Likelihood Estimator(MLE) of the parameter $\alpha$ we have to estimate it together with the scale and the location parameters. The corresponding MLE system of equations has no closed form solution. Prescott and Walden (1980)\cite{prescott1980maximum} proposed the last approach under more general settings, and more precisely for the family of all Generalized Extreme Values(GEV) distributions. The Probability of Waited Moments system of equations is derived by Hosking et al. (1985) \cite{hosking1985estimation}. The authors apply numerical methods in order to obtain its solution. In 2015 de Haan and Ferreira \cite{ferreira2015block}, suppose that the c.d.f. of the observed r.v.  belongs to the max-domain of attraction of some GEV distribution and investigate the Block maxima estimator for the shape parameter.  The corresponding MLE asymptotic theory  was recently developed by Dombry and Ferreira \cite{dombry2019maximum}. Due to the generality of their assumptions, however, their method is applicable only in cases of huge samples. B{\"u}cher and Segers (2018) \cite{bucher} generalize the results to stationary time series data, with distribution which belongs to the max-domain of attraction of Fr$\acute{e}$chet distribution. They describe the disadvantages of Peaks-over-threshold and Block maxima methods. The authors point out that these methods "...arise from asymptotic theory and are not necessarily accurate at sub-asymptotic thresholds or at finite
block lengths." Here we propose a simple parametric estimator $Q_{i,s}^{Fr*}$, defined in (\ref{QstarFr}) for the parameter $\alpha$. It is invariant with respect to a scale change of the sample.  The proof of the following theorem could be seen in Jordanova (2020) \cite{Jordanova2020Monograph}.

{\bf Theorem 4.}\label{thmTh1Frechetcase} Let us consider a sample of $n$ independent observations of a r.v. $X \in$ Fr$\acute{e}$chet($\alpha$, 0, 1), $\alpha > 0$.
  \begin{description}
 \item[i)] For any $s = 2, 3, ...$, $i = 2, 3, ...$,  and $y > 0$,
$$f_{\alpha Q_{i,s}^{Fr*}}(y) = c_s C_{i,s}e^{-yc_s} \int_{0}^\infty z\exp(-iz)\left\{\exp[-ze^{-yc_s}]-\exp(-z)\right\}^{i(s-1)-1}\left\{1 - \exp[-ze^{-yc_s}]\right\}^{i-1} \exp[-ze^{-yc_s}]dz,$$
where $c_s = -\log\left[1-\frac{\log(s)}{\log(s+1)}\right]$ and $C_{i,s} = \frac{[i(s+1)-1]!}{[(i-1)!]^2[i(s-1)-1]!}$. If $y \leq 0$, $f_{\alpha Q_{i,s}^{Fr*}}(y) = 0$.
 \item[ii)] $Q_{i,s}^{Fr*}$ estimator is strongly consistent, i.e. $Q_{i,s}^{Fr*}  {\mathop{\to}\limits_{i \to \infty}^{a.s.}} \frac{1}{\alpha}$.
 \end{description}

 In order to estimate quantiles outside the range of the data in this case we use the following estimator
 \begin{equation}\label{Froutsidetherange}
 \left(\hat{F}_{i,s,Fr}^\leftarrow\right)^*(p) = [-\log(p)]^{-Q_{i,s}^{Fr*}}.
 \end{equation}

 The proof of the following theorem is analogous to the one of Theorem 3.

 {\bf Theorem 5.}\label{QuantileEstimatorsFrpropertiesFr} Let us fix $s = 2, 3, ...$ and consider a sample of $n = i(s+1)-1$, $i \in \mathbb{N}$ independent observations on a r.v. $X \in$ Fr$\acute{e}$chet($\alpha$, 0, 1), $\alpha > 0$. Then
 $\left(\hat{F}_{i,s,Fr}^\leftarrow\right)^*(p)  {\mathop{\to}\limits_{i \to \infty}^{a.s.}} F_X^\leftarrow\left(p\right)$.

{\it Simulation study}

\begin{figure}
\begin{minipage}[t]{0.5\linewidth}
    \includegraphics[scale=.58]{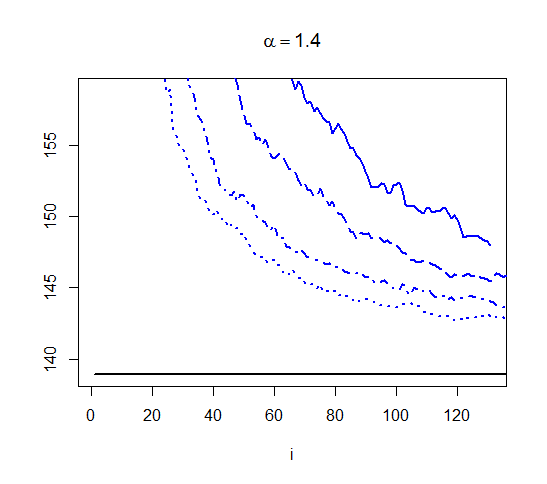}\vspace{-0.3cm}
\end{minipage}
\begin{minipage}[t]{0.49\linewidth}
    \includegraphics[scale=.58]{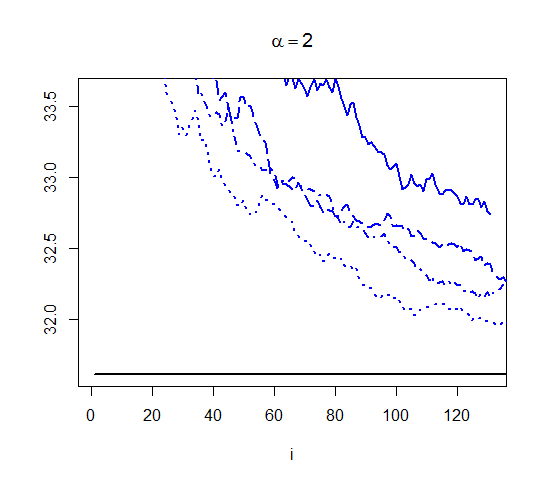}\vspace{-0.3cm}
\end{minipage}
\caption{\small Dependence of $\left(\hat{F}_{i,s,Fr}^\leftarrow\right)^*(p)$ on $i$ in Fr$\acute{e}$chet case.\label{fig:FrAMiTaNSquantiles}}
\end{figure}

The limit behaviour of $Q_{i,s}^{Fr*}$ is investigated and partially depicted in Jordanova (2020) \cite{Jordanova2020Monograph}. Along the current study we observed that for $\alpha < 1$ and $p = \frac{999}{1000}$ given small samples again $\left(\hat{F}_{i,s,Fr}^\leftarrow\right)^*(p)$ estimators are applicable only for very large samples. Analogously to the previous case we have plotted the images in Figure \ref{fig:FrAMiTaNSquantiles} for $\alpha = 1.4$ or $\alpha = 2$. The estimators $\left(\hat{F}_{i,s,Fr}^\leftarrow\right)^*(p)$ for $s = 2$ are plotted by solid lines. The case $s = 3$ is depicted by dashed lines. For $s =  4$ we have used dash-dot lines. The last case  $s =  5$ is presented via dotted lines.  The straight solid line visualise the estimated value of $F_X^\leftarrow\left(p\right)$. We observe again the strong consistency of the estimators and can conclude that they have very similar properties of the corresponding estimators considered in the Log-Logistic case in the previous section.

\section{HILL-HORROR CASE}
Let $\alpha > 0$. Embrechts et al. (2013) \cite{EKM} define the following distribution via its quantile function, $F_X^\leftarrow(p) = \frac{-\log(1 - p)}{(1 - p)^{1/\alpha}}$, $p \in (0, 1)$. Briefly we will denote it by $X \in HH(\alpha)$.  This distribution is called Hill-horror distribution with parameter $\alpha$, because of the difficulties related with the estimation of its parameter $\alpha$ via Hill\cite{hill1975simple} estimator. In this case, the distribution of $\alpha \log(X)$ depends on $\alpha$. Therefore, if we consider  logarithms of ratios of two order statistics and compute their means (if they exist) the result will depend non-linearly on $\alpha$. In order to obtain a strongly consistent estimator of $\frac{1}{\alpha}$, before we normalise the logarithm of fractions of order statistics with $\log(s)$, we have shifted its distribution with an appropriate constant. It is determined by the equality
$$\log\left[\frac{F_X^\leftarrow\left(\frac{s}{s+1}\right)}{F_X^\leftarrow\left(\frac{1}{s+1}\right)}\right] = \frac{1}{\alpha}\log(s) - \log\left[1 - \frac{\log(s)}{\log(s + 1)}\right].$$
In this way we obtain $Q_{i,s}^{HH*}$ estimator, defined by formula (\ref{QstarHH}). When we would like to estimate the quantiles outside the range of the data we obtain $\left(\hat{F}_{i,s,HH}^\leftarrow\right)^*(p) = \frac{-\log(1 - p)}{(1 - p)^{Q_{i,s}^{HH*}}}$.
Now, the strong consistency of these estimators follows by a Continuous mappings theorem.

{\bf Theorem 6.}\label{thmTh1HHcase}  Let us fix $s = 2, 3, ...$ and consider a sample of $n = i(s+1)-1$, $i \in \mathbb{N}$ independent observations on a r.v. $X \in HH(\alpha)$, $\alpha > 0$. Then,
\begin{description}
\item[i)] $Q_{i,s}^{HH*}  {\mathop{\to}\limits_{i \to \infty}^{a.s.}} \frac{1}{\alpha}$;
\item[ii)] $\left(\hat{F}_{i,s,HH}^\leftarrow\right)^*(p)  {\mathop{\to}\limits_{i \to \infty}^{a.s.}} F_X^\leftarrow\left(p\right)$.
\end{description}

\begin{figure}
\begin{minipage}[t]{0.5\linewidth}
    \includegraphics[scale=.58]{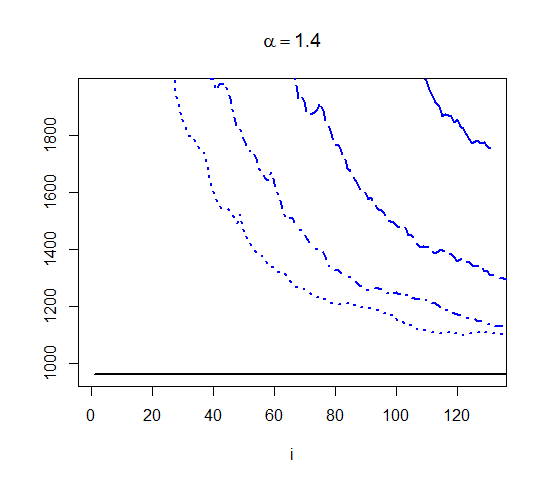}\vspace{-0.3cm}
\end{minipage}
\begin{minipage}[t]{0.49\linewidth}
    \includegraphics[scale=.58]{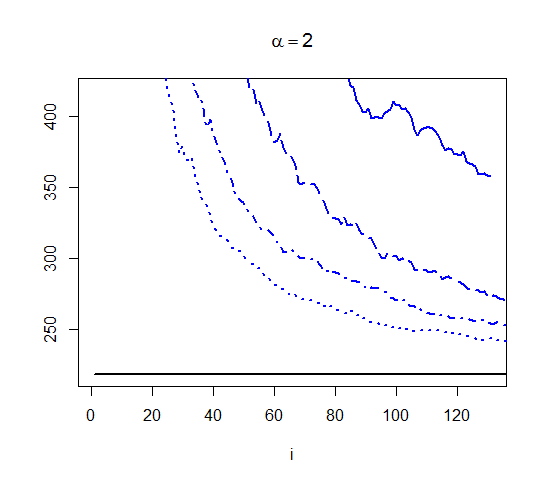}\vspace{-0.3cm}
\end{minipage}
\caption{\small Dependence of $\left(\hat{F}_{i,s,HH}^\leftarrow\right)^*(p)$ on $i$ in Hill-horror case.\label{fig:HHAMiTaNSquantiles}}
\end{figure}

{\it Simulation study}
Jordanova (2020) \cite{Jordanova2020Monograph} compares this estimator with Hill\cite{hill1975simple} estimators and shows that $Q_{i,s}^{HH*}$ estimator has a relatively fast rate of convergence. In this section we will see that although $\bar{F}_X \in RV_{-\alpha}$ and the statistic $\left(\hat{F}_{i,s,HH}^\leftarrow\right)^*(p)$ is again a strongly consistent estimator of $F_X^\leftarrow\left(p\right)$, its rate of convergence to the quantiles outside the range of the data is very slow. Our simulation study shows very different results than in the previous two cases. These considerations speak once again that although we have fixed $\alpha$, the class of c.d.f. with regularly varying right tails with this parameter $\alpha$ is too wide in order to be possible the parameter $\alpha$ to be simultaneously non-parametrically estimated within this class. In order to plot our parametric estimators $\left(\hat{F}_{i,s,HH}^\leftarrow\right)^*(p)$ we have followed a similar procedure for drawing the graphs that was used in the previous two sections. We have simulated 1000 samples on $n = 150(s+1)-1$ independent observations on $X \in HH(\alpha)$ separately for $\alpha = 1.4$ and $\alpha = 2$. For any fixed $\alpha$, $p = \frac{999}{1000}$ and $s = 2, 3, 4, 5$, and  for any $i = 1, 2, ..., 150$ we have computed $\left(\hat{F}_{i,s,HH}^\leftarrow\right)^*(p)$ estimators. Then, for any fixed $s$ and $i$ and over these $1000$ samples we have determined the average of $\left(\hat{F}_{i,s,HH}^\leftarrow\right)^*(p)$ estimators. Finally, we have plotted them in Figure \ref{fig:HHAMiTaNSquantiles}. Again we observe the strong consistency of the considered estimators, however the rate of convergence is much slower than in previous two cases. This is very important especially when $\alpha$ approaches $0$. The last means that these last estimators are appropriate only if we are working with large samples.
The sample size $n = 150(s+1)-1$, where $s = 2, 3, 4, 5$, is not enough to observe their good properties.

An analogous approach could be applied to many different cases of distributions with regularly varying right tails. It leads us to strongly consistent and distribution sensitive estimators of the index of regular variation. These could be for example Exponentiated-Fr$\acute{e}$chet, Burr, Reverse Burr, Danielson and de Vries, among others distributions. If the observed r.v. does not have c.d.f. with regularly varying tail, however it can be continuously transformed to such a distribution the same approach is applicable to the transformed r.vs. These are for example $H_1$, Log-Pareto, or Weibull distributions.

\section{ACKNOWLEDGMENTS}
The authors are grateful to the bilateral projects Bulgaria - Austria, 2016-2019, Feasible statistical modelling for extremes in ecology and finance, Contract number 01/8, 23/08/2017 and WTZ Project BG 09/2017.

\end{document}